\documentclass[12pt,onecolumn]{IEEEtran}

\usepackage{multirow}
\usepackage{latexsym}
\usepackage{amsmath}
\usepackage{amsfonts}
\usepackage{amssymb}
\usepackage{graphicx}
\usepackage{times}
\usepackage{cite}
\usepackage{grffile}
\usepackage{amsmath}
\usepackage{times}
\usepackage{float}

\begin{document}

$$\textbf{Signal Detection in Singular Value Decomposition}$$

$$\text{Mohsen Rakhshan}$$

\section{Abstract}
We develop an Iterative version of the Singular Value Decomposition (ISVD) that jointly analyzes a finite number of data matrices to identify signals that correlate among the rows of matrices. It will be illustrated how the supervised analysis of a big data set by another complex, multi-dimensional phenotype using the ISVD algorithm could lead to signal detection.

\section{Introduction}
The joint analysis of different types of big data should provide a better understanding of the molecular mechanisms underlying cancers \cite{az1}, \cite{az2}, \cite{az3}, \cite{Huang}, \cite{Alter} and \cite{Stranger}. In this work we introduce an approach that involves the use of matrix factorizations such as the Singular Value Decomposition (SVD) to simultaneously analyze a finite number of matrices by iteratively identifying cross-correlations between them. The Iterative Singular Value Decomposition (ISVD) jointly analyzes a finite number of data matrices by automatically identifying any correlations that may exist among the rows of matrices, or the rows and columns of any one of the matrices. Since big data sets contain many distinct signals that associate the measured variables and samples, repeated application of the ISVD on re-normalized data iteratively detects these signals as orthogonal correlations among the matrices or within a single matrix. This approach is computationally efficient. The following Notation is used in this work.

\textbf{Notation:}
\begin{itemize}
\item $R^m:$ The linear space of real column vectors with $m$ coordinates.
\item $M_{mn}(R):$ The linear space of real $m \times n$ matrices.
\item $O_{m}(R):$ The subset of $m \times m$ real orthonormal matrices.
\item $D_{mn}(R):$ The subset of $m \times n$ real non-negative diagonal matrices.
\item $N_{mn}(R):$ The subset of $m \times n$ real matrices whose $n$ columns are normalized.
\end{itemize}

\section{Singular Value Decomposition (SVD)}
The singular value decomposition of a matrix $A \in M_{mn}(R)$ is the decomposition of $A$ into the product of three matrices of the form
$$A = UDV^T$$
where $U \in O_{m}(R)$, $D \in D_{mn}(R)$, and $V \in O_{n}(R)$. The diagonal entries of $D$ are called the singular values of $A$, and columns of $U$ and $V$ are called left and right singular vectors of $A$ respectively. The left-singular vectors are eigenvectors of $AA^T$ and right-singular vectors are eigenvectors of $A^TA$. Non-zero singular values of $A$ are same as the square roots of non-zero eigenvalues of $AA^T$ (which are the same as the non-zero eigenvalues of $A^TA$). Based on the SVD, the matrix $A$ can be written as an outer product of orthogonal rank one approximations
$$A=\sum_{i=1}^{\text{rank($A$)}} (u_iv_i^T) d_i  $$
where $u_i$ and $v_i$ are the $i$th columns of $U$ and $V$ respectively; $d_i$ is the $i$th diagonal element of $D$.

\subsection{Signal detection in SVD}
Given a set of $N$ matrices $A_i \in M_{m_i n}$ with full column rank. The SVD of stacked matrix $\begin{bmatrix} A_1 \\ \vdots \\ A_N \end{bmatrix}$ will be given by
\begin{equation}\label{eq:svd}
\begin{bmatrix} A_1 \\ \vdots \\ A_N \end{bmatrix}=P D Q^T, P \in O_{(m_1+\cdots+m_N)}(R), D \in D_{(m_1+\cdots+m_N)n}(R), Q \in O_{n}(R).
\end{equation}
The $j$th column of $P$ can be written as
$$p_j=\begin{bmatrix} p_{1j} \\ \vdots \\ p_{Nj} \end{bmatrix}$$
where $p_{ij} \in R^{m_i}$. If $q^j$ represent the $j$th ($j \leq n$) columns of $Q$, then $p_{1j},...,p_{Nj}$ can be viewed as correlation detectors of the signal $q^j$ that is common to $A_1,...,A_N$. That is, we can detect the rows of $A_1,...,A_N$ that are mutually correlated with the common signal $q^j$ by applying thresholds to the vectors $p_{1j},...,p_{Nj}$. Note that the above procedure will result in $n$ distinct signals that link the rows of $A_1,...,A_N$ via the common signals $q^j$. The analysis of the rows of $A_1,...,A_N$ will be based on $Q$.

\section{Iterative Singular Value Decomposition Algorithm}\label{se:isvd}
Specifically, given the rank one approximation to the data we form the residual data and apply rank one approximation to the residual data pair. This procedure is iterated until the largest singular value of residue falls below a certain threshold that is near zero. This algorithm is based on the SVD of the stacked matrix $\begin{bmatrix} A_1 \\ \vdots \\ A_N \end{bmatrix}$. Define $q_i$ ($p_i$) to be the $i$th column of matrix $Q$ ($P$) and $d_i$ be the $i$th diagonal entry of matrix $D$ in (\ref{eq:svd}). The first column of $P$ can be written as
$$p_1=\begin{bmatrix} p_{11} \\ \vdots \\ p_{N1} \end{bmatrix}$$
where $p_{i1} \in R^{m_i}$. The rank one approximation to SVD of $\begin{bmatrix} A_1 \\ \vdots \\ A_N \end{bmatrix}$ is given as
$$\begin{bmatrix} A_1^1 \\ \vdots \\ A_N^1 \end{bmatrix}=\begin{bmatrix} p_{11}  d_1  q_1^T \\ \vdots \\ p_{N1}  d_1 q_1^T \end{bmatrix}=\begin{bmatrix} p_{11} \\ \vdots \\ p_{N1} \end{bmatrix} d_1  q_1^T=p_1 d_1  q_1^T,$$
$$A_1^1 = p_{11}/|p_{11}| . |p_{11}| . d_1q_1^T     = u_{11} . d_{11} . v_1^T,$$
$$\vdots$$
$$A_N^1 = p_{N1}/|p_{N1}| . |p_{N1}| . d_1q_1^T = u_{N1} . d_{N1} . v_1^T,$$
with
$$d_{11}^2+... + d_{N1}^2 = 1,$$
$$v_1^T v_1 = d_1^2.$$
We are interested to find the support of signals in $A_1,...,A_N$. The vectors $u_{11}=\frac{p_{11}}{\| p_{11} \|},...,u_{N1}=\frac{p_{N1}}{\| p_{N1} \|}$ in the rank one approximation of $A_1,...,A_N$ can be viewed as detectors of the signal in $A_1,...,A_N$, respectively. That is, the components of $u_{11},...,u_{N1}$ with large absolute magnitude correspond to the rows in $A_1,...,A_N$ that are highly correlated with $v_1^T=d_1  q_1^T$. Note that $v_1$ is common in the matrix factorization of $A_1,...,A_N$. The indices of $u_{11},...,u_{N1}$ with large absolute magnitude are the signal support for $A_1,...,A_N$ respectively. We apply the universal threshold proposed by Donoho and Johnstone (\cite{Donoho04} and \cite{Donoho05}) to detect components of $u_{11},...,u_{N1}$ with large absolute magnitude.

Now we subtract $A_1^1,...,A_N^1$ from $A_1,...,A_N$ respectively and get the rank one approximation to SVD of the pair $\begin{bmatrix} A_1-A_1^1 \\ \vdots \\ A_N-A_N^1 \end{bmatrix}$. It is shown in Theorem~$1$ that
$$\begin{bmatrix} A_1-A_1^1 \\ \vdots \\ A_N-A_N^1 \end{bmatrix}=\sum_{i=2}^{n} p_i d_i  q_i^T.$$
The rank one approximation to SVD of $\begin{bmatrix} A_1-A_1^1 \\ \vdots \\ A_N-A_N^1 \end{bmatrix}$ is given as
$$\begin{bmatrix} A_1^2 \\ \vdots \\ A_N^2 \end{bmatrix}=\begin{bmatrix} p_{12}  d_2  q_2^T \\ \vdots \\ p_{N2}  d_2 q_2^T \end{bmatrix}=\begin{bmatrix} p_{12} \\ \vdots \\ p_{N2} \end{bmatrix} d_2  q_2^T=p_2 d_2  q_2^T,$$
$$A_1^2 = p_{12}/|p_{12}| . |p_{12}| . d_2q_2^T = u_{12} . d_{12} . v_2^T,$$
$$\vdots$$
$$A_N^2 = p_{N2}/|p_{N2}| . |p_{N2}| . d_2q_2^T = u_{N2} . d_{N2} . v_2^T,$$
with
$$d_{12}^2+... + d_{N2}^2 = 1,$$
$$v_2^T v_2 = d_1^2.$$
$$v_2^T v_1 = 0.$$
Similarly to the first iteration, the signals in $u_{12}=\frac{p_{12}}{\| p_{12} \|},...,u_{N2}=\frac{p_{N2}}{\| p_{N2} \|}$ are detected. The sequential application of the ISVD is continued until the maximum singular value of
$$\begin{bmatrix} A_1-A_1^1-A_1^2-\cdots-A_1^k \\ \vdots \\ A_N-A_N^1-A_N^2-\cdots-A_N^k \end{bmatrix} $$
drops below a predefined threshold indicating that the noise floor of the data has been reached. The number of iterations is at most rank of $\begin{bmatrix} A_1 \\ \vdots \\ A_N \end{bmatrix}$. Therefore, the ISVD algorithm is computationally efficient and enhances signal detection by rank one approximation. To get the largest singular value and singular vector of $A$ we can use the MATLAB command $eig(A^TA,1)$. For further information about the algorithm of finding a few eigenvalue of a large matrix please refer to \cite{Lehoucq}.

\textbf{Theorem 1:}

Let the SVD of the stack matrix be
\[  \begin{bmatrix} A_1  \\ \vdots \\ A_N \end{bmatrix} = PDQ^T = \sum_{i=1}^r p_i d_i q_i^T . \]
Then, splitting apart the first column of $P$, we can write
\begin{align*}
A_1^1 = p_{11} d_1 q_1^T &= \frac{p_{11}}{\| p_{11} \|} \cdot \| p_{11} \| \cdot d_1 q_1^T = u_{11} \cdot d_{11} \cdot v_1^T \\
   \vdots\\
A_N^1 = p_{N1} d_1 q_1^T &= \frac{p_{N1}}{\| p_{N1} \|} \cdot \| p_{N1} \| \cdot d_1 q_1^T = u_{N1} \cdot d_{N1} \cdot v_1^T
\end{align*}
satisfying
\begin{align*}
&d_{11}^2+ \cdots + d_{N1}^2 = 1 \\
&v_1^T v_1 = d_1^2 \\
&u_{11}^T u_{11} = \cdots =u_{N1}^T u_{N1} = 1 .
\end{align*}
Let
\[ H = \begin{bmatrix} A-A_1^1 \\ \vdots \\ A-A_N^1 \end{bmatrix} = \sum_{i=2}^r p_i d_i q_i^T . \]
which latter equality holds since
\[ \begin{bmatrix} A_1^1  \\ \vdots \\ A_N^1 \end{bmatrix} = p_1 d_1 q_1^T . \]
The claim is that
\begin{itemize}
\item the nonzero eigenvalues of $HH^T$ are $p_2, \dots, p_r$,
\item the nonzero eigenvalues of $H^TH$ are $q_2, \dots, q_r$,
\item with the corresponding eigenvalues $d_2^2, \dots ,d_r^2$,
\end{itemize}
so that the above representation of $H$ is its singular value decomposition.

\textbf{Proof:} The proof is in the calculation
\begin{align*}
HH^Tp_2 = (\sum_{i=2}^r p_id_iq_i^T)(\sum_{j=2}^r q_jd_jp_j^T) p_2 = (\sum_{i=2}^r p_id_iq_i^T)q_2 d_2 = d_2^2 p_2
\end{align*}
and similarly $H^THq_2 = d_2^2 q_2$, etc.

The algorithm proceeds to find
\begin{align*}
&A_1^2 = p_{12} d_2 q_2^T = \frac{p_{12}}{\| p_{12} \|} \cdot \| p_{12} \| \cdot d_2 q_2^T = u_{12} \cdot d_{12} \cdot v_2^T& \\
&\vdots&\\
&A_N^2 = p_{N2} d_2 q_2^T = \frac{p_{N2}}{\| p_{N2} \|} \cdot \| p_{N2} \| \cdot d_2 q_2^T = u_{N2} \cdot d_{N2} \cdot v_2^T&
\end{align*}
satisfying
\begin{align*}
&d_{12}^2 + \cdots+ d_{N2}^2 = 1 \\
&v_2^T v_2 = d_2^2 \\
&u_{12}^T u_{12}=\cdots=u_{N2}^T u_{N2} = 1 \\
&v_1^T v_2 = 0 .
\end{align*}

\textbf{Theorem 2:} Let $s_i$ be an $n \times 1$ unit vector and $\{s_i\}_{i=1}^{2}$ be a set of orthonormal vectors. Define two matrices $A_1$ ($m_1 \times n$) and $A_2$ ($m_2 \times n$) as follows
\[ A_1=\begin{bmatrix} s_1^T  \\ \vdots \\ s_1^T \\ s_2^T  \\ \vdots \\ s_2^T  \end{bmatrix}, \]
and
\[ A_2=\begin{bmatrix} s_1^T  \\ \vdots \\ s_1^T \end{bmatrix}, \]
where the support of $s_1$ and $s_2$ in $A$ are $m_{s_1}$ and $m_{s_2}$ respectively and the support of $s_1$ in $A_2$ is $m_2$. Define $H$ as the stack of two matrices $A_1$ and $A_2$
\[ H_1=\begin{bmatrix} A_1  \\ A_2  \end{bmatrix}, \]
It is claimed that
\begin{itemize}
\item $m_{s_1}+m_{2}$ and $m_{s_2}$ are the eigenvalues of $H_1^T H_1$ corresponding to eigenvector $s_1$ and $s_2$ respectively.
\item If $m_{s_1}+m_{2}=m_{s_2}$ then $s_1+s_2$ and $s_1-s_2$ are the eigenvectors of $H_1^T H_1$ corresponding to eigenvalue $m_{s_1}+m_{2}=m_{s_2}$.
\end{itemize}
Proof:
\begin{equation}\label{eq1}
H_1^T H_1=A_1^TA_1+A_2^TA_2=\begin{bmatrix} s_1  & \cdots & s_1 &s_2 & \cdots & s_2  \end{bmatrix}\begin{bmatrix} s_1^T  \\ \vdots \\ s_1^T \\ s_2^T  \\ \vdots \\ s_2^T  \end{bmatrix}+\begin{bmatrix} s_1 & \cdots & s_1 \end{bmatrix}\begin{bmatrix} s_1^T  \\ \vdots \\ s_1^T \end{bmatrix}
\end{equation}
$$=(m_{s_1}+m_2) s_1 s_1^T+m_{s_2}s_2 s_2^T$$
By multiplying $s_1$ to both sides of (\ref{eq1})
$$H_1^T H_1s_1=(m_{s_1}+m_2) s_1 s_1^Ts_1+m_{s_2}s_2 s_2^Ts_1=(m_{s_1}+m_2) s_1$$
Similarly for $s_2$.

Suppose $m_{s_1}+m_{2}=m_{s_2}$, by multiplying $s_1$ to both sides of (\ref{eq1})
$$H_1^T H_1(s_1+s_2)=(m_{s_2} s_1 s_1^T+m_{s_2}s_2 s_2^T)(s_1+s_2)=m_{s_2} s_1 s_1^Ts_1+m_{s_2}s_2 s_2^Ts_1+m_{s_2} s_1 s_1^Ts_2+m_{s_2}s_2 s_2^Ts_2$$
$$=m_{s_2} s_1+m_{s_2}s_2=m_{s_2} (s_1+s_2).$$
Similarly for $s_1-s_2$.

\textbf{Theorem 3:} Let $s$ be an $n \times 1$ unit vector and $\{n_i\}_{i=1}^{p}$ be a set of i.i.d. white noise vectors, in which $n_i$ is a $n \times 1$ vector. Define the matrices $A$ ($p \times n$) as follows
\[ A=\begin{bmatrix} s^T+n_1^T  \\ \vdots \\ s^T+n_m^T \\ n_{m+1}^T\\ \vdots \\ n_{p}^T \end{bmatrix}, \]
It is claimed that
\begin{itemize}
\item $m$ is the approximate eigenvalues of $A^T A$ corresponding to eigenvector $s$ for a non-trivial range of signal to noise ratio.
\end{itemize}
Proof:
\begin{flalign} \label{eq:eq1111}
A^TA &=\begin{bmatrix} s+n_1  & \cdots & s+n_m & n_{m+1} & \cdots & n_{p}\end{bmatrix}\begin{bmatrix} s^T+n_1^T \\ \vdots \\ s^T+n_m^T \\ n_{m+1}^T \\ \vdots \\  n_{p}^T \end{bmatrix} \\ \notag
&=\,\ \,\ ss^T+n_1s^T+sn_1^T+n_1n_1^T\\ \notag
& \,\ \,\ \,\ \,\ \,\ \,\ \,\ \,\ \,\ \,\ \,\ \,\ \,\ \,\ \,\ \vdots \\ \notag
&\,\ \,\ \,\ +ss^T+n_ms^T+sn_m^T+n_mn_m^T \\ \notag
&\,\ \,\ \,\ +n_{m+1}n_{m+1}^T+\cdots+n_{p}n_{p}^T
\end{flalign}
By multiplying $\frac{s}{m}$ from right hand side to both sides of (\ref{eq:eq1111})
\begin{flalign} \label{eq:eq2222}
A^TA\frac{s}{m}=s+\frac{\sum_{i=1}^{m}n_i}{m}
\end{flalign}
It is easy to see that if
\begin{flalign}
\|s\| \gg \|\frac{\sum_{i=1}^{m}n_i}{m}\|
\end{flalign}
or
\begin{flalign}\label{eq:eq22221}
\frac{<s,s>}{\sigma_n^2} \gg \frac{n}{m}
\end{flalign}
then the right hand side of (\ref{eq:eq2222}) can be approximated by $s$. It is noticeable that if $m \gg n$, then with a high probability the inequality in (\ref{eq:eq22221}) is satisfied. Below the convergence behaviour of $A^TA\frac{s}{m}$ as $m$ grows is studied.

\subsection{convergence behaviour}
By taking limit from both sides of (\ref{eq:eq2222}) we have
\begin{flalign} \label{eq:eq222222}
\lim A^TA\frac{s}{m}=s+\lim \frac{\sum_{i=1}^{m}n_i}{m}=s
\end{flalign}
By law of large number we have
\begin{flalign} \label{eq:eq3333}
Pr(\lim \frac{\sum_{i=1}^{m}n_i}{m}=0)=1
\end{flalign}
From (\ref{eq:eq222222}) and (\ref{eq:eq3333}) we can have almost surely
\begin{flalign} \label{eq:eq4444}
\lim A^TA\frac{s}{m}=s
\end{flalign}
Therefore almost surely for any given $\epsilon>0$ there exists an $M$ such that for all $m>M$
\begin{flalign} \label{eq:eq5555}
|A^TA\frac{s}{m}-s|<\epsilon
\end{flalign}
and almost surely
\begin{flalign} \label{eq:eq6666}
m(s-\epsilon)< A^TAs< m(s+\epsilon).
\end{flalign}
From Chebychev's inequality
\begin{flalign} \label{eq:eq6666}
P(|\frac{1}{m}\sum_{i=1}^{m}n_i|>\epsilon) \leq \frac{var(\frac{1}{m}\sum_{i=1}^{m}n_i)}{\epsilon^2}=\frac{\sigma_n^2}{m\epsilon^2}
\end{flalign}
Let $\epsilon$ be the weakest signal to noise ratio that thresholding strategy works, $\epsilon=\frac{m}{\|N\|}$.
\begin{flalign} \label{eq:eq6666}
P\big(|\frac{1}{m}\sum_{i=1}^{m}n_i|>\frac{m}{\|N\|}\big) \leq \frac{\sigma_n^2\|N\|^2}{m^3}=\frac{\sigma_n^4pn}{m^3}=\frac{\sigma_n^4pn}{m^3}
\end{flalign}

What is the critical $\epsilon$? Critical $\epsilon$ can be achieved from ROC curve and Donoho Johnson threshold. For critical $\epsilon$ we can find an upper bound on probability of error for any value of error. The maximum probability of error is a function of signal to noise ratio, because the critical $\epsilon$ is a function of signal to noise ratio and $m$ is the power of signal.

\section{Conclusion}

The ISVD algorithm is able to detect signals that are common or unique in finite number of data sets. Signals are detected sequentially based on signal strength. Very weak signals with small singular values are detectable since noise background is systematically reduced relative to signal strength. Note that the ISVD is computationally feasible in situations where the GSVD is not. It may also be theoretically more meaningful than the Generalized Singular Value Decomposition (GSVD) \cite{VanLoan}, \cite{Page}, \cite{Friedland05} and \cite{Orly}. 

\end{document}